\documentclass[10pt,a4paper]{article}
\usepackage[left=24mm,right=24mm,top=20mm,bottom=20mm]{geometry}
\usepackage{cmap}
\usepackage[T1]{fontenc}
\usepackage{lmodern}
\usepackage[english]{babel}
\usepackage{amsmath,amssymb,amsthm,booktabs,microtype}
\usepackage{titlesec}
\titlespacing*{\section}{0pt}{2.4ex plus .5ex minus .2ex}{1.25ex plus .15ex}
\usepackage[hidelinks]{hyperref}
\newcommand{\F}{\mathbb{F}}
\newcommand{\eps}{\varepsilon}
\newtheorem{theorem}{Theorem}
\newtheorem{lemma}[theorem]{Lemma}
\newtheorem{corollary}[theorem]{Corollary}
\title{Strongly regular graphs on $1296$ vertices\\[2pt]
and two infinite families}
\author{Mikhail P. Golubyatnikov}
\date{}

\begin{document}
\maketitle
\vspace{-2.9em}
\begin{center}
\small
\textsuperscript{1}Krasovskii Institute of Mathematics and Mechanics UB RAS,\\
16 S.~Kovalevskaya Str., Yekaterinburg, 620077, Russia\\
\textsuperscript{2}Ural Mathematical Center, 16 S.~Kovalevskaya Str.,
Yekaterinburg, 620077, Russia\\
\textsuperscript{3}Ural Federal University named after the First President of
Russia B.~N.~Yeltsin, Yekaterinburg, Russia\\[3pt]
Email: \texttt{mike\_ru1@mail.ru}
\end{center}
\vspace{0.15em}

\begin{abstract}
We construct new strongly regular graphs from graphs with parameters of
negative Latin square type. The construction yields one graph
with parameters $(1296,407,118,132)$ and another with parameters
$(1296,444,144,156)$. Each family contains at least $9!=362880$ pairwise non-isomorphic graphs with automorphism group $C_3$. Iteration gives two infinite families, each containing
a graph of order $1296\cdot9^j$ for every $j\geq0$.
\end{abstract}



\section{Introduction}

In Brouwer's table, as consulted on 28 November 2026, both parameter
sets, $(1296,407,118,132)$ and $(1296,444,144,156)$, are listed
as open~\cite{table}. We give a common construction for these two parameter
sets. In both cases, three graphs on a common $144$-element set $X$
determine adjacency between vertices with equal first coordinates, and
adjacency between vertices with distinct first coordinates is given by
the same function $f\colon X\to\F_3$. Once the parameters of
the three graphs are fixed, the conditions in Theorem~\ref{thm:criterion}
depend only on $f$. Successive Wang--Qiu--Hu switches on pairs chosen from
nine fixed triples give at least $9!$ isomorphism classes for each parameter set
(Lemma~\ref{lem:isomorphisms}).

All graphs are finite, simple, and undirected. A graph on $v$ vertices
is \emph{strongly regular with parameters $(v,k,\lambda,\mu)$} if it
is regular of valency $k$, where $0<k<v-1$, and two distinct vertices
have exactly $\lambda$ common neighbours when adjacent and exactly
$\mu$ when nonadjacent. We write $N_G(u)$ for the neighbourhood
of $u$ in $G$ and $\F_p$ for the field of $p$ elements.

\section{The construction}

Let $X$ be a finite abelian group of order $n^2$, written additively.
Fix a function $f\colon X\to\F_3$, three graphs $B_0,B_1,B_2$ with
vertex set $X$, and $\eps\in\{0,1\}$. The vertex set of $\Gamma$
is $\F_3\times X\times\F_3$. Two distinct vertices are adjacent
according to the following rule.
\[
 (i,x,s)\sim(j,y,t)\quad\Longleftrightarrow\quad
 \begin{cases}
   t-s=f(y-x), & j=i+1,\\
   s-t=f(x-y), & j=i-1,\\
   x\sim_{B_i}y\ \text{or}\ (x=y\ \text{and}\ \eps=1), & j=i .
 \end{cases}
\]
Write $V_i=\{i\}\times X\times\F_3$ and $D_{i,x}=\{i\}\times\{x\}\times\F_3$.
Each $D_{i,x}$ induces a triangle when $\eps=1$ and an independent
set when $\eps=0$. The edges between $D_{i,x}$ and $D_{i+1,y}$
form the perfect matching given by $s\mapsto s+f(y-x)$.

\section{A criterion for strong regularity}

\begin{theorem}\label{thm:criterion}
Suppose that $3\mid n$, let $a=n/3-\eps$, and assume that each
$B_i$ is strongly regular with parameters
\[
 (n^2,b,\lambda_0,\mu_0)
   =\bigl(n^2,\,a(n+1),\,a^2+3a-n,\,a(a+1)\bigr).
\]
Then $\Gamma$ is strongly regular if and only if both of the
following conditions hold.
\begin{enumerate}
\item\label{cond:diff} For every $d\in X\setminus\{0\}$, the
  difference $f(z+d)-f(z)$ takes each value in $\F_3$ exactly
  $n^2/3$ times as $z$ runs through $X$.
\item\label{cond:sum} For every $d\in X$, the expression
  $f(z)+f(d-z)+f(-d)$ takes the values $0,1,2$ with respective
  multiplicities $(n^2-2n)/3,\ (n^2+n)/3,\ (n^2+n)/3$ as $z$
  runs through $X$.
\end{enumerate}
Under these conditions, the parameters of $\Gamma$ are
\begin{equation}\label{parameters}
 \bigl(9n^2,\ (n-\eps)(3n+1),\ (n-\eps)^2-3\eps,\
              (n-\eps)(n-\eps+1)\bigr).
\end{equation}
\end{theorem}

\begin{proof}
Every vertex of $V_i$ has $3b+2\eps$ neighbours in $V_i$ and $n^2$
in each of $V_{i-1}$ and $V_{i+1}$. Suppose first that $\Gamma$ is
strongly regular, with common-neighbour numbers $\lambda$ and $\mu$.
The vertices $(i,x,s)$ and $(i,x,t)$, with $s\ne t$, have
$3b+\eps$ common neighbours, all in $V_i$. Hence $\mu=3b$
if $\eps=0$, and $\lambda=3b+1$ if $\eps=1$.

For a pair with distinct first coordinates, write the vertices as $(i,x,s)$ and
$(i+2,y,t)$. Put $d=y-x$ and $h=t-s+f(-d)$. The pair is
adjacent precisely when $h=0$. For each neighbour $w$ of $x$ in
$B_i$, the pair has one common neighbour in $D_{i,w}$.
The unique neighbour of $(i+2,y,t)$ in $D_{i,x}$ differs from
$(i,x,s)$ precisely when $h\ne0$. It is then a common neighbour
if $\eps=1$. The same count in $V_{i+2}$ gives a total of
$2b$ in $V_i\cup V_{i+2}$ when $h=0$ and $2b+2\eps$
when $h\ne0$.

A vertex $(i+1,x+z,r)$ is a common neighbour precisely when
\[
 r=s+f(z)=t-f(d-z).
\]
The number of common neighbours in $V_{i+1}$ is therefore the
number of $z\in X$ satisfying $f(z)+f(d-z)+f(-d)=h$.
For fixed $d$, summing over $h\in\F_3$ counts each $z$ once
and gives $\lambda+2\mu=6b+4\eps+n^2$. Together with the value of
$\mu$ or $\lambda$ obtained above, this identity yields
\eqref{parameters}. The numbers of solutions for $h=0$ and for each
$h\ne0$ are $\lambda-2b$ and $\mu-2b-2\eps$, respectively.
Substitution gives condition~\ref{cond:sum}.

Finally, consider two vertices $(i,x,s)$ and $(i,y,t)$ with
$d=y-x\ne0$. In $V_i$, they have $3\lambda_0+4\eps$
common neighbours if $x\sim_{B_i}y$, and $3\mu_0$ otherwise.
In each of $V_{i-1}$ and $V_{i+1}$ the count is $D_d(t-s)$, where
\[
 D_d(h)=\bigl|\{z\in X:f(z+d)-f(z)=h\}\bigr|.
\]
Indeed, for a common neighbour with second coordinate $w$, take
$z=w-y$ in $V_{i+1}$ and $z=x-w$ in $V_{i-1}$.
Since $\lambda-3\lambda_0-4\eps=\mu-3\mu_0=2n^2/3$, we obtain
$D_d(h)=n^2/3$ for every $h$. This proves
condition~\ref{cond:diff}.

Conversely, suppose both conditions hold and take $\lambda,\mu$
from~\eqref{parameters}. The counts for the three types of pairs above
give $\lambda$ common neighbours for every adjacent pair and $\mu$
for every nonadjacent pair. The valency computed at the outset is
$(n-\eps)(3n+1)$.
\end{proof}

\section{Graphs on 1296 vertices}

Take $n=12$, $T=\F_2^2\times\F_3$, and $X=T^2$. Set
$f(x_1,x_2)=p(x_1)+p(x_2)$, where $p\colon T\to\F_3$.
The entry in row $(a,b)$ and column $c$ gives $p(a,b,c)$.
\[
 \begin{array}{c|ccc}
 (a,b)\backslash c & 0 & 1 & 2\\ \hline
 (0,0) & 2 & 0 & 2\\
 (0,1) & 0 & 1 & 2\\
 (1,0) & 0 & 0 & 0\\
 (1,1) & 2 & 1 & 0
 \end{array}
\]
The twelve values of $p$ give the following multiplicities.
\begin{center}
\begin{tabular}{@{}llrrr@{}}
\toprule
Expression, $z\in T$ & Condition on $u$ & $0$ & $1$ & $2$\\
\midrule
$p(z+u)-p(z)$       & $u\ne0$ & 4 & 4 & 4\\
$p(z)+p(u-z)+p(-u)$ & $u\in T$ & 4 & 2 & 6\\
\bottomrule
\end{tabular}
\end{center}
For $d=(d_1,d_2)\ne0$, at least one coordinate is nonzero. The
first row then gives each value of $f(z+d)-f(z)$ exactly
$4\cdot12=48$ times. For condition~\ref{cond:sum}, we add two
residues with multiplicities $(4,2,6)$ in each coordinate.
The residue $0$ occurs $4^2+2\cdot2\cdot6=40$ times. The multiplicities
of $1$ and $2$ are $2\cdot4\cdot2+6^2=52$ and $2\cdot4\cdot6+2^2=52$, respectively.
These are the required values for $n=12$.

For each value of $\eps$, take $B_0,B_1,B_2$ to be copies of a
single graph $B$. For $\eps=1$, we use the graph of valency $39$
whose edges form an orbit of $\mathrm{PGL}_3(3)$ acting on unordered
pairs of right cosets of a subgroup of order $39$~\cite{orbital}.
For $\eps=0$, we use Brady's Cayley graph of valency
$52$~\cite{brady}. Theorem~\ref{thm:criterion} gives the following
parameters.
\begin{center}
\begin{tabular}{@{}crll@{}}
\toprule
$\eps$ & $a$ & Parameters of $B$ & Parameters of $\Gamma$\\
\midrule
1 & 3 & $(144,39,6,12)$  & $(1296,407,118,132)$\\
0 & 4 & $(144,52,16,20)$ & $(1296,444,144,156)$\\
\bottomrule
\end{tabular}
\end{center}

For each $B$, the labelling of $X$ is the one recorded in the accompanying
graph6 files~\cite{catalogue}, with the triples $(i,x,s)$ taken in
lexicographic order.

\section{Pairwise non-isomorphic graphs}
\label{sec:nonisomorphism}

Put $0_T=(0,0,0)$, $e=(0,1,0)$ and $\tau(x)=x+(e,e)$. Let
\[
 U=\{0_T\}\times\{(a,b,c)\in T:(a,b)\ne(1,1)\}\subset X .
\]
Extend each permutation $\pi$ of these nine points by the identity on
$X\setminus U$. Write $\Gamma_\pi$ for the graph obtained with
$B_0=\pi(B)$ and $B_1=B_2=\tau(B)$, where $\pi(B)$ denotes the graph
obtained by relabelling $x$ as $\pi(x)$. The map $x\mapsto(1,\tau(x),0)$
identifies the labelled graph $B$ with an induced subgraph of every $\Gamma_\pi$.
For fixed $B$, the induced graphs $\Gamma_\pi[V_0]$ are isomorphic.
The matchings between distinct $V_i$ are prescribed by the same $f$ for every $\pi$.

\begin{lemma}\label{lem:isomorphisms}
For each of the parameter sets $(1296,407,118,132)$ and $(1296,444,144,156)$,
there are at least $9!=362\,880$ pairwise non-isomorphic strongly regular
graphs with full automorphism group isomorphic to $C_3$. Each is obtained
from the corresponding $\Gamma_{\mathrm{id}}$ by a sequence of Wang--Qiu--Hu switches.
\end{lemma}

\begin{proof}
Let $A$ be the adjacency matrix of $\pi(B)$ and take $P=D_{0,x}$,
$Q=D_{0,y}$ for distinct $x,y\in U$. The subgraphs induced by $P$, $Q$
and $P\cup Q$ are regular of degrees $2\eps$, $2\eps$ and $2\eps+3A_{xy}$.
For $v\notin P\cup Q$,
\[
 \bigl(|N_{\Gamma_\pi}(v)\cap P|,\,|N_{\Gamma_\pi}(v)\cap Q|\bigr)=
 \begin{cases}
 (1,1), & v\notin V_0,\\
 3(A_{zx},A_{zy}), & v=(0,z,s).
 \end{cases}
\]
Wang--Qiu--Hu switching~\cite{wqh} interchanges the neighbour counts
$(3,0)$ and $(0,3)$ and preserves all other adjacencies~\cite[Theorem~1.2]{ihringer}.
This transposes $x$ and $y$ in $B_0$ and takes $\Gamma_\pi$ to $\Gamma_{(x\,y)\circ\pi}$.
Transpositions generate $\operatorname{Sym}(U)$, and
Theorem~\ref{thm:criterion} gives the parameters of every $\Gamma_\pi$.

Write $\mathcal D=\{D_{i,x}:i\in\F_3,\ x\in X\}$.
For $a=1,2$, let $H_a$ be the spanning subgraph of $\Gamma_\pi$ obtained by
deleting the edges within each $V_i$ with $i\ne a$. It is independent of $\pi$.

Let $S=\{v_1,v_2,v_3\}\subseteq V(\Gamma_\pi)$ be any set of three distinct
vertices. Write $N_G(S)=\bigcap_{v\in S}N_G(v)$ for their common neighbourhood in $G$.
If $S\notin\mathcal D$, put $I(S)=\{i\in\{0,2\}:S\cap V_i\ne\varnothing\}$.
For each $i\in I(S)$ choose $u\in S\cap V_i$ and $v\in S$ in a different
triple of $\mathcal D$. By the pair counts in the proof of Theorem~\ref{thm:criterion},
$u,v$ have at most $L=\max\{b+\eps,3\lambda_0+4\eps,3\mu_0\}$ common neighbours in $V_i$.
Hence
\[
 |N_{\Gamma_\pi}(S)|
 =|N_{H_1}(S)|+\sum_{i\in I(S)}|N_{\Gamma_\pi}(S)\cap V_i|
 \le |N_{H_1}(S)|+|I(S)|L .
\]
For $m=0,1,2$, let $M_m$ be the largest value of $|N_{H_1}(S)|$ over
triples $S\notin\mathcal D$ with $|I(S)|=m$. For each $S$, the value of
$|N_{H_1}(S)|$ was found by counting the columns of the adjacency matrix
of $H_1$ that contain a $1$ in all three corresponding rows.
\[
 \begin{array}{c|rrr|r|rr}
 k & M_0 & M_1 & M_2 & L & \max_m(M_m+mL) & 3b \\ \hline
 407 & 87 & 72 & 28 & 40 & 112 & 117\\
 444 &105 & 72 & 33 & 60 & 153 & 156
 \end{array}
\]
The common neighbours of $D_{i,x}$ are the vertices of the $b$ triples
$D_{i,y}$ with $y\in N_{B_i}(x)$. The last two columns therefore give,
for every triple $S$,
\[
 S\in\mathcal D\ \iff\ |N_{\Gamma_\pi}(S)|=3b .
\]
Two distinct triples $D_{i,x}$ and $D_{j,y}$ have exactly three edges
between them if and only if $i\ne j$. Thus every isomorphism permutes the sets $V_i$.

Put $\rho(i,x,s)=(i,x,s+1)$ and $R(i,x,s)=(i+1,x,s)$.
The permutation $\rho$ preserves $\Gamma_\pi$ and $H_1$.
For each $k$, \texttt{nauty}~2.9.3 gives $|\operatorname{Aut}(H_1)|=3$
with no prescribed vertex colours~\cite{nauty}. Since $H_2=R(H_1)$
and $R\rho=\rho R$, both $H_1$ and $H_2$ have automorphism group $\langle\rho\rangle$.

Let $g:\Gamma_\pi\to\Gamma_\sigma$ be an isomorphism. At most one of $V_1,V_2$
maps to $V_0$, so $g(V_a)=V_d$ for some $a,d\in\{1,2\}$. Since
$g(H_a)=H_d=R^{d-a}(H_a)$, the map $R^{a-d}g$ belongs to $\langle\rho\rangle$.
Thus for some $c\in\F_3$,
\[
 g(i,x,s)=(i+d-a,x,s+c).
\]
Take $p=((0,0,1),0_T)$ and $q=(e,0_T)$, both outside $U$. In both cases,
the edge counts between $D_{i,p}$ and $D_{i,q}$ form the vector $(0,9,9)$
in the order $i=0,1,2$. The unique zero forces $d=a$, so $g=\rho^c$
and $\pi(B)=\sigma(B)$. For $x\ne y$,
\[
 |N_B(x)\mathbin{\triangle}N_B(y)|
 =2\bigl(b-|N_B(x)\cap N_B(y)|\bigr)
 \ge 2\bigl(b-\max\{\lambda_0,\mu_0\}\bigr)>|U|.
\]
Thus $x\mapsto N_B(x)\setminus U$ is injective. The automorphism
$\sigma^{-1}\pi$ of $B$ fixes $X\setminus U$ pointwise, so $\pi=\sigma$.
We obtain
\[
 \operatorname{Iso}(\Gamma_\pi,\Gamma_\sigma)=
 \begin{cases}
 \langle\rho\rangle,&\pi=\sigma,\\
 \varnothing,&\pi\ne\sigma.
 \end{cases}
\]
The $9!$ permutations of $U$ give pairwise non-isomorphic graphs,
each with full automorphism group $\langle\rho\rangle\cong C_3$.
\end{proof}

\begin{corollary}\label{cor:automorphisms}
For each $\pi$, the group $\operatorname{Aut}(\Gamma_\pi)$ is generated by
$(i,x,s)\mapsto(i,x,s+1)$. Its $432$ vertex orbits are the triples $D_{i,x}$.
\end{corollary}

\section{Iteration and two infinite families}

Suppose $X,f,B_0,B_1,B_2$ satisfy the conditions of Theorem~\ref{thm:criterion}.
Let $X'=X\times\F_3^2$ and identify $(i,x,s)\in V(\Gamma)$
with $(x,i,s)\in X'$. For the next step take $B'_0=B'_1=B'_2=\Gamma$ and set
\[
 f'(x,u,v)=f(x)+uv .
\]
Put $n'=3n$. Then $|X'|=(n')^2$ and $a'=n'/3-\eps=n-\eps$.
In particular, $a'^2+3a'-n'=a'^2-3\eps$, so the parameters
in~\eqref{parameters} are exactly those required of each $B'_i$.

Write $Q(u,v)=uv$ and fix $(d,p)\in X'\setminus\{0\}$.
If $p=0$, then $d\ne0$, and each value of $f(z+d)-f(z)$ is
repeated nine times. If $p\ne0$, write
$w=(w_1,w_2)$. The additional term in the difference is
\[
 Q(w+p)-Q(w)=p_2w_1+p_1w_2+p_1p_2 .
\]
Its linear part is nonzero. For each fixed $z\in X$, the difference
therefore takes each value three times as $w$ runs through $\F_3^2$.
In either case the multiplicity of each value is $(n')^2/3$.

For any $(d,p)\in X'$, the additional term in condition~\ref{cond:sum} is
\[
 Q(w)+Q(p-w)+Q(-p)=-Q(w+p)
\]
by expansion in $\F_3$. As $w$ varies, this term takes the value
$0$ five times and each nonzero value twice. If $c_h$ is the
multiplicity of $h$ in the second condition for $f$, its
multiplicity for $f'$ is consequently
\[
 5c_h+2(n^2-c_h)=2n^2+3c_h .
\]
This gives $((n')^2-2n')/3$ for $h=0$ and $((n')^2+n')/3$
for $h=1,2$. Thus both conditions persist under iteration.

Starting with $n=12$ and keeping $\eps$ fixed, we obtain the
parameters in~\eqref{parameters} with $n=12\cdot3^j$ for every
$j\geq0$. The two families have order $1296\cdot9^j$ at step $j$.
For $j=1$, their parameters are $(11664,3815,1222,1260)$ for
$\eps=1$ and $(11664,3924,1296,1332)$ for $\eps=0$.
For $j\geq1$, the equality $B'_0=B'_1=B'_2$
makes $R$ an automorphism as well, and $\langle R,\rho\rangle\cong C_3\times C_3$.
The least eigenvalue for the parameters in~\eqref{parameters} is
$\theta=-(2n+\eps)$.

\end{document}